\documentclass{amsart}

\usepackage[T1]{fontenc}
\usepackage[latin1]{inputenc}

\usepackage{amsmath}
\usepackage{amssymb}
\usepackage{amsfonts}

\makeatletter

\usepackage{graphicx}

\makeatother

\begin{document}

\title{Gaussian Summation:\\An exponentially convergent summation scheme}
\author{
  H. Monien
}
\address{
  Physikalisches Insitut\\
  Universit\"at Bonn\\
  Nu{\ss}allee 12\\
  53115 Bonn\\
  Germany
}
\email{monien@th.physik.uni-bonn.de}
%\homepage{http://cond-mat.uni-bonn.de}
\date{\today}
\thanks{}
%\affiliation{Physikalisches Institut, Universit\"at Bonn}
\begin{abstract}
  Gaussian Quadrature is a well known technique for numerical
  integration.  Recently Gaussian quadrature with respect to discrete
  measures corresponding to finite sums have found some new interest.
  In this paper we apply these ideas to infinite sums in general and
  give an explicit construction for the weights and abscissae of {\em
    Gaussian summation} formulas. The abscissae of the Gaussian
  summation have a very interesting asymptotic distribution function
  with a (cusp) singularity.  We apply the Gaussian summation
  technique to two problems which have been discussed in the
  literature. We find that the Gaussian summation has an extremely
  rapid convergence rate for the Hardy-Littlewood sum for a large
  range of parameters.  For functions which are smooth but have a
  large scale $a$ the error of Gaussian Summation shows
  exponential convergence as a function of summation points, $n$: the
  error decreases like $\sim\; n\exp(-4n^2/\pi a)$.  The Gaussian
  summation achieves a given accuracy with $n\sim\sqrt{a}$ points
  whereas other summation schemes require at least $n\sim a$ function
  evaluations.
\end{abstract}
\maketitle

\section{introduction}
Scalar sums of the form
\begin{equation}
  \label{eq:sum}
  \sum_{n=-\infty}^{\infty}f(n) 
\end{equation}
with $n$ integer and $f$ being a continuous function appear in nearly every
context of mathematics and physics. In many problems, either because $f$ is
for example the solution to another complicated nonlinear problem or there
exists no known analytical expression for the sum, one has to resort to
numerical summation techniques.  Because of its importance many techniques
have been developed for a fast evaluation of Eq.~\ref{eq:sum} which depend on
the asymptotic behavior the function $f(x)$. A good summary of the classical
techniques can be found in \cite{BenderOrzag}. Newer developments include
application of Levin-type convergence transformation schemes (\cite{Levin} for
a review see \cite{Homeier}).  In the application of these schemes the
sequence of partial sums $s_n=\sum_{m=1}^{n} f(m)$ is calculated and then
extrapolated to $n\rightarrow\infty$. However if the function $f$ has a large
intrinsic scale so that the asymptotic behavior is reached only at very large
values of $n$ then these schemes are inefficient because a large number of
terms has to be accumulated before the asymptotic behavior is reached and not
applicable because rounding errors prevent a reliable extrapolation.  So the
problem of evaluating the sum Eq.~\ref{eq:sum} efficiently and reliably for
problems with a large intrinsic scale is interesting.

In this paper we discuss a summation scheme based on ideas related to Gauss
Quadrature.  The basic idea is to replace the sum Eq.~\ref{eq:sum} by another
sum
\begin{equation}
  \label{eq:sum_gauss}
  \sum_{n=-\infty}^{\infty}f(n) \approx \sum_{k=1}^{N} w_k f(x_k)
\end{equation}
where the weights $w_k$ and abscissae $x_k$ are chosen in such a way
as to approximate closely the sum Eq.~\ref{eq:sum} for a large class
of functions for a relatively small number $N$. Because the values of
$x_k$ are not fixed to be integer numbers the hope is that the
asymptotic behavior at large $x$ is captured by the sum
Eq.~\ref{eq:sum_gauss}. The central question is what are the optimal
weights, $w_k$ and points $x_k$?  This is the question we are going to
address in this paper. The plan of the paper is as follows: First we
derive the moment generating Weyl function for sums of the form
$\sum_{\nu\in\Omega} 1/\nu^2 f(1/\nu^2)$ and use the Pade approximants
of the Weyl functions to obtain the orthogonal polynomials related to
these sums.  For numerical application we give the Jacobi matrix for
generating the weights and moments for Gaussian summation formulas and
give two examples of the application of the Gaussian summation.
Finally we derive the asymptotic distribution of the zeros of the
above orthogonal polynomials.

Similar schemes have been reviewed recently by Engblom \cite{Engblom}.
Milovanovi{\'c} and Cvectkovi{\'c} \cite{Milovanovic} studied the convergence
properties of Gaussian summation schemes for sums of the form
$\sum_{\nu\ge0}p^{-\nu}f\left(a q^{-\nu}\right)$, with $|p|, |q|>1$.  The
earliest reference to Gaussian summation can be found in the book by Uvarov
and Nikiforov as an example using Chebyshev polynomials for finite sums.

%--- Moment generating functions ----------------------------------------------

\section{moment generating functions and recursion relation}

We define a linear functional $L[p]$, acting on the space of
polynomials $\mathcal{P}$, using some discrete measure $\lambda$. 
\begin{equation}
  \label{eq:general_functional}
  L[p] = \int p d\lambda 
\end{equation}
The focus of the paper is on sums of the type 
\begin{equation}
  \label{eq:functional}
  L[f] = \sum_{\nu\in\Omega}
  \frac{1}{\nu^2}\;f\left(\frac{1}{\nu^2}\right)
\end{equation}
where $\Omega$ is a subset of $\mathbb{Z}$. 
The moments are defined as 
\begin{equation}
  \label{eq:moments}
  \mu_n = \int x^n d\lambda = \sum_{\nu\in\Omega} \frac{1}{\nu^{2(n+1)}} 
\end{equation}
The moment generating Weyl function for Eq.~4 is given by
\begin{equation}
  \label{eq:Weyl_function}
  \Phi(z) = \int\frac{d\mu}{z-x} = 
  \sum_{\nu\in\Omega} \frac{1}{\nu^2}\frac{1}{z-\frac{1}{\nu^2}} = 
  \frac{1}{z}\sum_{k=0}^\infty \frac{\mu_k}{z^k}.
\end{equation}
Weyl function can be calculated explicitly for
$\Omega=\mathbb{Z}\backslash\{0\}$:
\begin{equation}
  \label{eq:weyl_function_all}
  \Phi(z) =  1-\frac{\pi}{\sqrt{z}}\cot\left(\frac{\pi}{\sqrt{z}}\right),
\end{equation}
Eq.~\ref{eq:weyl_function_all} serves as starting point for the
deriving the orthogonal polynomials for the corresponding weight
function.
%--- Recursion relation -------------------------------------------------------
The Pade approximant of the moment generating function are closely
related to the orthogonal polynomials of the corresponding weight
function \cite{Baker}.  Here we derive an explicit expression for the
Pade approximant of Eqs.~\ref{eq:weyl_function_all} which is
equivalent to determining the orthogonal polynomials for the weight
function \cite{Baker}.  Using the continued fraction expansion of
$\cot(x)$ \cite{Abramowitz} we can express the generating function as:
\begin{equation}
  \label{eq:continued_fraction_cot}
  \Phi(x=\pi/\sqrt{z}) = 1-x\cot(x)=
  \frac{x^2}{3}
  \cfrac{1}{1-
    \cfrac{x^2/15}{1-
      \cfrac{x^2/35}{1-
        \cfrac{x^2/63}{1-
          \ldots
        }
      }
    }
  }
\end{equation}
The numerator of the $n$th Pade approximant corresponding to these
continued fractions, $P_n(x)$, is denoted by $R_n(x)$ and the
denominator by $S_n(x)$.  The denominator and numerator obey the
recursion relations \cite{BenderOrzag}:
\begin{eqnarray}
  \label{eq:recursion_pade}
  R_{n+1}(x)&=&R_{n}(x)+c_{n+1}\;x^2\;R_{n-1}(x)\\
  S_{n+1}(x)&=&S_{n}(x)+c_{n+1}\;x^2\;S_{n-1}(x)
\end{eqnarray}
for $n=0, 1, 2 \ldots$ with the initial condition for the recursion
are $R_{-1}(x)=0$, $R_0(x)=c_0 x^2$, $S_{-1}(x)=1$ and $S_0(x)=1$.  The
coefficients $c_n$ are given by $c_0=1/3$,
$c_n=-1/(2n+1)(2n+3)$, $n=1,2,\ldots$. 

Surprisingly enough the recursion relation for $R_n(x)$ and $S_n(x)$
can be solved analytically. Let us discuss the recursion for Eq.
\ref{eq:continued_fraction_cot}.  First note that the recursion
for the Bessel functions, $Z_\nu$ of half-integer index $n+1/2$ are
given by:
\begin{equation}
  \label{eq:Bessel_recursion}
  Z_{n+3/2}(x)-\frac{2n+1}{x}Z_{n+1/2}(x)+Z_{n-1/2}(x)=0.
\end{equation}
Splitting of the asymptotic behavior for large $n$ we define $a_n(x)$
by requiring that
\begin{equation}
  \label{eq:Bessel_asymptotic}
  Z_{n+1/2}(x) = 
  \left(\frac{2}{x}\right)^n\Gamma(n+\frac{1}{2})
  \;a_n(x).
\end{equation}
The coefficients $a_n(x)$ obey the recursion relation:
\begin{equation}
  \label{eq:recursion_a_n}
  a_{n+1}=a_{n}-\frac{x^2}{4n^2-1}\;a_{n-1}
\end{equation}
This is precisely the recursion relation for the $S_n(x)$ if we shift
the index $n$ by 2. Thus the recursion relation for $S_n(x)$ has a
solution of the form:
\begin{equation}
  \label{eq:S_solution_general}
  S_n(x) = \frac{x^{n+2}}{2^{n+2}\Gamma(n+\frac{5}{2})}
  \left[
    A(x)J_{n+5/2}(x)+B(x)Y_{n+5/2}(x),
  \right]
\end{equation}
where $A(x)$ and $B(x)$ are two unknown functions which only depend on $x$ and
can be determined from the initial conditions. The orthogonal polynomials with
respect to the functional Eq.  \ref{eq:functional} are given by denominators
of the $[M-1/M]$ Pade approximant \cite{Baker} which are $\left\{S_{-1}(x),
  S_1(x), S_3(x), S_5(x), \ldots\right\}$. Transforming back to the original
variable $z = (\pi/x)^2$ we define the polynomials
\begin{equation}
  \label{eq:orthogonal_polynomial}
  s_n(z) = z^n S_{2n-1}\left(\frac{\pi}{\sqrt{z}}\right).
\end{equation}
These are the sought after orthogonal polynomials for the functional,
Eq.~\ref{eq:functional}. Specializing to the initial conditions an
explicit solution for $S_n(x)$ can be obtained:
\begin{equation}
  \label{eq:solution_S}
  S_n(x=\pi/\sqrt{z})=
  -\frac{\pi x^{n+3/2}}{\Gamma(n+\frac{5}{2})2^{n+5/2}}
    \left[
     \cos(x)J_{n+5/2}(x)+\sin(x)Y_{n+5/2}(x).
   \right]
\end{equation}
This result will be useful to determine the distribution of zeros explicitly.
The first few polynomials for $n=0, 1, \ldots$ are given by
\begin{eqnarray*}
  \label{eq:first_S}
  S_0(x) &=& 1\\
  S_1(x) &=& 1-\frac{1}{15}x^2\\
  S_2(x) &=& 1-\frac{2}{21}x^2\\
  S_3(x) &=& 1-\frac{1}{9}x^2+\frac{1}{945}x^4\\
  S_4(x) &=& 1-\frac{4}{33}x^2+\frac{1}{495}x^4\nonumber\\
  &\ldots&
\end{eqnarray*}
We note that using Eq.~\ref{eq:solution_S} it can be shown easily that the
polynomials $s_n(x)$ are related to Bessel polynomials, $y_n(x)$, of imaginary
argument \cite{Nikiforov}:
\begin{equation}
  \label{eq:bessel_polynomials}
  s_n(x) = \frac{\sqrt{\pi}}{2\;4^n\Gamma(2n+\frac{3}{2})}
  {\rm Re}\left(y_{2n+1}\left(-i\frac{\sqrt{x}}{\pi}\right)\right).
\end{equation}
Below we will derive simple three-point recursion relations for the 
polynomials $s_n(x)$.

%--- calculation of the weights and abscissae

\section{Calculation of the weights and abscissae}
The orthogonal polynomials with respect to the functional
Eq.~\ref{eq:general_functional} are $\{S_0, S_1, S_3, S_5, \ldots \}$.  It is
useful to obtain a recursion relation for the odd numbered polynomials only.
This can be achieved by eliminating the even numbered polynomials from the
recursion.  After some algebra we find:
\begin{equation}
  \label{eq:s_recursion_2}
  S_{n+3}(x)=\left(1-\frac{2x^2}{(2n+5)(2n+9)}\right)S_{n+1}(x) -
  \frac{x^4}{(2n+3)(2n+5)^2(2n+7)}S_{n-1}(x).
\end{equation}
This can be brought into the usual three-term recursion relation,
${s}_{k+1}(x) = (x-a_k){s}_k(x) - b_k{s}_{k-1}(x)$. The coefficients
of the recursion for $k=1, 2, \ldots$ are given by:
\begin{eqnarray}
  \label{eq:recursion_coefficients}
  a_k &=& \frac{2\pi^2}{(4k+1)(4k+5)}\\
  b_k &=& \frac{\pi^4}{(4k-1)(4k+1)^2(4k+3)}.
\end{eqnarray}
The first coefficients have to be determined from the initial conditions and
are $a_0=\pi^2/15$ and $b_0$ is arbitrary.  With these coefficients the problem
of finding the weights and abscissae of the Gaussian summation formulas
reduces to the standard procedure used to determine the analogous quantities
in Gaussian integration \cite{Golub}. The Jacobi matrix for determining the
weights and abscissae is given by
\begin{equation}
  \label{eq:Jacobi_matrix}
  \tilde J_N=\left[
    \begin{array}{ccccccc}
      a_0 & \sqrt{b_1} &&&&& \\
      \sqrt{b_1} & a_1 & \sqrt{b_2} &&&& \\
      & \sqrt{b_2} & a_2 & \sqrt{b_3} &&& \\
      &&\ddots&\ddots&\ddots&&\\
      &&&\sqrt{b_{N-2}}&a_{N-2}&\sqrt{b_{N-1}}\\
      &&&&\sqrt{b_{N-1}}&a_{N-1}
    \end{array}
  \right]
\end{equation}
It is easy to calculate the weights and abscissae using an implementation of
the Golub and Welsh algorithm \cite{Golub}, either ``gauss.m'' from the QCP
collection of Gautschi \cite{Gautschi}, or the {\em gaucof} routine from the
{\em Numerical Recipes} \cite{NumericalRecipes}. In terms of the weights
$w_n$ and abscissae $x_n$ determined from the Jacobi matrix,
Eq.~\ref{eq:Jacobi_matrix} the Gaussian summation procedure is 
given by: 
\begin{equation}
  \label{eq:Gaussian_summation}
  \sum_{\nu\in\Omega} \frac{1}{\nu^2} f\left(\frac{1}{\nu^2}\right) =
  \sum_{k=1}^n {w}_k f({x}_k) + \epsilon_n
\end{equation}
with the error estimate $\epsilon_n$ (see e.g. \cite{Engblom}). The 

\begin{equation}
  \label{eq:error}
  \epsilon_n = K_n\frac{f^{(2n)}(\xi)}{(2n)!}
\end{equation}
where $\xi\in C(\Omega)$ and $C(\Omega)$ is defined as the smallest
connected set containing $\Omega$. The constant $K_n$ can be
calculated using Eq. (2.17) in \cite{Engblom}:
\begin{equation}
  \label{eq:K_n}
  K_n = ||P_n||^2 = \frac{1}{2} (4n+3)\pi^{4n+3}
    16^{-(n+1)}\;\Gamma^{-2}\left(2n+\frac{5}{2}\right)  
\end{equation}
Although the error estimate Eq.~\ref{eq:error} gives an upper bound it
is not very useful to estimate the convergence of the Gaussian
summation formulas. For the specific examples below we can sharpen the
error estimates and achieve a more general understanding of the
convergence of the Gaussian summation formulas. 
%---- section numerical results -----------------------------------------------
\section{Two examples}
We apply the result of the preceeding section to two examples. 
The first is the Hardy-Littlewood function which has interesting
properties in itself and is defined as:
\begin{equation}
  \label{eq:hardy_littlewood_sum}
  H(x)=\sum_{k=1}^\infty \frac{\sin(x/k)}{k}.
\end{equation}
This summand has an expansion in terms of $1/k^2$ for $k/x\gg 1$. For 
large values of $x$ the summation is not expected to converge before
the $k\approx x$. The numerical evaluation of this function has been
discussed by Gautschi \cite{Gautschi,GautschiHL}. The methods considered
in these references are Gaussian integration of the Laplace transformation
of the sum and the application of Euler-MacLaurin asymptotic summation
the eliminate the first asymptotic terms $O(1/k^2)$ and $O(1/k^4)$.

\begin{figure}
  \centering
  \includegraphics[angle=-90,width=\textwidth]{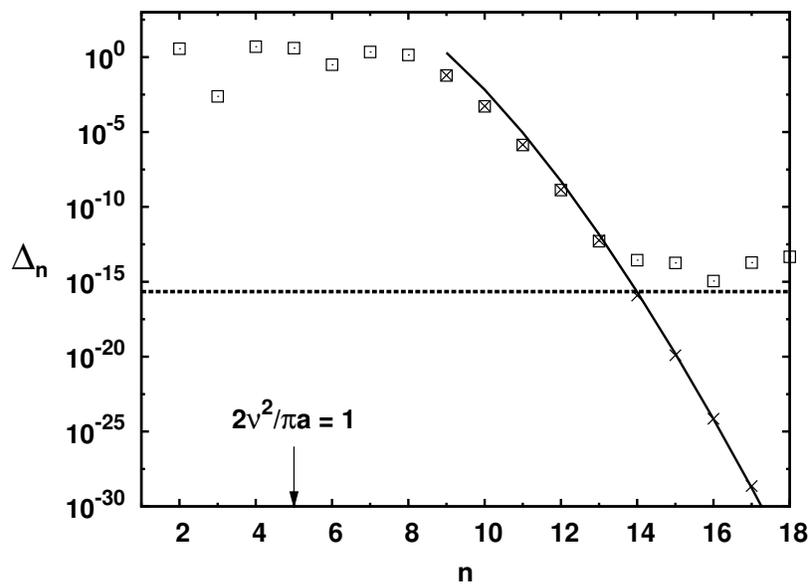}
  \caption{Relative error in the evaluation of the Hardy-Littlewood sum
    for $x=100$. The squares ($\Box$) are indicating the numerical
    error in the Gaussian summation, the crosses ($\times$) are
    indicating the saddle point approximation of the error integral
    and the solid line gives the asymptotic behavior of the former for
    large $n$.  The dashed line gives the machine~$\epsilon$ for
    comparison.}
  \label{fig:hardy_littlewood_error}
\end{figure}

We have applied our method to problem Eq.~\ref{eq:hardy_littlewood_sum} for
$x=100$. The summation converges rapidly as shown below and for $n>14$ the
error in the sum is dominated by the error in the weights and points. One
observes that the differences of the value of the Gaussian summation for
different $n$ is fluctuating around the machine $\epsilon$. We find that the
Gaussian summation scheme achieves a high accuracy using very few function
evaluation of the summand. Basically the accuracy is limited by the accuracy
of the Gaussian weights and points used in the evaluation of the summation.
For example using fifteen terms the Gaussian summation gives the result for
the Hardy-Littlewood sum with an error less then $10^{-14}$ for all
$x\in(0,100)$.  To compare specifically to the results by Gautschi in
reference \cite{Gautschi}, Table 3.21, with the largest value of $x$, $x=40$,
we compare the number of evaluations needed to achieve a relative accuracy of
$\epsilon_0=2.22\times 10^{-7}$ as stated there.  According to this reference
the smallest number of points in the integration of the Laplace transformed
sum for this accuracy is $N=39$ so a total of $39\times70$ (for the
integration) function evaluations are needed.  With the Gaussian summation
method a total $n=8$ function evaluations are needed to achieve the same
relative accuracy for $H(40)$. We summarize our results in
Table~\ref{tab:hardy_littlewood_sum_error}.  The symbol $-$ in the table
indicates when the relative error in the summation was less then $10^{-14}$.
The error is plotted in Fig.~\ref{fig:hardy_littlewood_error} for $x=100$ as a
function of $n$. Gaussian summation give results already at a relatively small
number of terms which does not grow like the scale where the function reaches
the asymptotic behavior.

Next we derive an analytic expression for the error estimate confirming the
empirical findings. Obviously
\begin{equation}
  \label{eq:contour_approx}
  \Sigma_n[f] = \sum_{k=1}^n w_k f(x_k) = 
  \frac{1}{2\pi i}\int_{C} d\zeta\;\frac{1}{\zeta}
  \left(1 - \frac{R_{2n-1}(\pi\zeta)}{S_{2n-1}(\pi\zeta)}\right)\;
  \frac{1}{\zeta^2}f\left(\frac{1}{\zeta^2}\right),
\end{equation}
where $C$ is a contour enclosing the zeros of $S_n(\pi\zeta)$.  The
error estimate $\Delta_n[f] = \Sigma_{n+1}[f]-\Sigma{n}[f]$ can be simplified
using the recursion Eq.~\ref{eq:recursion_pade}
\begin{equation}
  \label{eq:error_explicit}
  \Delta_n[f] = 
  \frac{1}{2\pi i} \int_{C} d\zeta\; A_n(\pi\zeta)\;
  \frac{1}{\zeta^2}f\left(\frac{1}{\zeta^2}\right),
\end{equation}
with 
\begin{equation}
  \label{eq:def_Q}
  A_n(\zeta) =  
  \frac{2\pi(4n+3)\pi^{4n+1}\zeta^{4n+1}}{2^{4n+5}\Gamma^2(2n+5/2)}
  \times
  \frac{1}{S_{2n+1}(\pi\zeta)S_{2n-1}(\pi\zeta)}
\end{equation}
The asymptotic behavior of the error is determined by
Eq.~\ref{eq:error_explicit} for $n\gg 1$. Using the asymptotic form of
the Bessel function on the imaginary axis \cite{Abramowitz} and
introducing $\nu=2n+5/2$ we obtain as the leading contribution:
\begin{equation}
  \label{eq:A_phi}
  A_n(i \nu\tau) \approx 8\;\nu e^{\Phi_\nu(\nu\tau)}
\end{equation}
with
\begin{equation}
  \label{eq:phi}
  \Phi_\nu(z) = \sum_\pm
  \left[
    -\frac{1}{2}\ln\left(\frac{z}{\nu\pm 1}\right)
    +\frac{1}{2}\ln\left(1+\left(\frac{z}{\nu\pm 1}\right)^2\right)
    -(\nu\pm 1)\eta\left(\frac{z}{\nu\pm 1}\right)
  \right]
\end{equation}
and
\begin{equation}
  \label{eq:eta}
  \eta(\tau)=\sqrt{\tau^2+1}+\log\left(\frac{\tau}{1+\sqrt{\tau^2+1}}\right).
\end{equation}
For the Hardy-Littlewood function the error can be estimated as:
\begin{equation}
  \label{eq:error_HL}
  \Delta_n \approx \frac{4\nu}{\pi}\int_C d\zeta\;
  \exp({\Phi_\nu(\zeta)})\; \sinh\left(\frac{\pi a}{\zeta}\right)
\end{equation}
This integral is completely dominated by the saddle points on the imaginary 
axis which appear for $2\nu^2/\pi a \gg 1$ and are located located at 
\begin{equation}
  \label{eq:tau_min}
  \nu\tau_{\rm min} = \pm\frac{\nu}{\sqrt{\xi^2-2\xi}}
  \left(
    1 + \frac{\xi}{2(\xi-2)\nu} + O\left(\frac{1}{\nu^2}\right)
  \right),
\end{equation}
with the direction of steepest descent being perpendicular to the
imaginary axis and $\xi$ is given by $\xi=2\nu^2/\pi a$. The 
evaluation of the integral at the saddle point yields:
\begin{equation}
  \label{eq:saddle_point_integral}
  \Delta_n\approx 
  8\nu\;
  \sqrt{\frac{2\pi}{\Phi_\nu^{\prime\prime}(\nu\tau_{\rm min})}}\;
  \exp\left(
    \Phi_\nu\left(\nu\tau_{\rm min}\right)
  \right)\;
  \sinh\left(
    \frac{\pi a}{\nu\tau_{\rm min}}
  \right).
\end{equation}
For $2\nu^2/\pi a \gg 1$ this expression can be simplified further to
finally yield the asymptotic error estimate:
\begin{equation}
  \label{eq:error_HL_asymptotic}
  \Delta_n \approx 2\sqrt{\pi\nu}
  \;
  e^{-\pi a/\nu}
  \;
  \left(
    \frac{e^2\pi a}{4\nu^2}
  \right)^{2\nu}.
\end{equation}
We have plotted this error estimate in Fig.~\ref{fig:hardy_littlewood_error}.
The agreement of the saddle point approximation
Eq.~\ref{eq:saddle_point_integral} and the analytic asymptotic expression
Eq.~\ref{eq:error_HL_asymptotic} with the numerical result is quite
reasonable. This expression confirms that for smooth functions the Gaussian
summation has a very rapid rate of convergence. The rate of convergence is
determined by the scale of the function - in this case $a$ at which the
asymptotic behavior sets in. However in the Gaussian summation the rate of
convergence is determined by the ratio $a/\nu^2$ so that convergence already
sets in at $\nu \approx \sqrt{a}$. This is a more rapid convergence rate than
any of the extrapolation schemes for sums proposed in \cite{GautschiHL}.

The second example serves to illustrate this point more clearly.  We
consider a sum known explicitly \cite{Abramowitz}:
\begin{equation}
  \label{eq:example_2}
  G(a)=\sum_{k=-\infty}^{\infty}\frac{1}{a^2+k^2} = \frac{\pi}{a}\coth(\pi a),
\end{equation}
with $a\in\mathbb{R}_{>0}$. We have applied Gaussian summation to
Eq.~\ref{eq:example_2}. Again we observe rapid convergence even for
values of $a$ as large as 1000. The numerical convergence pattern is
completely regular. In this case the expression for the error
estimator $\Delta_n$ can be given without approximation because the
contour integral can be done exactly by deforming the contour and is
given by the residues of the integrand at $\pm ia$:
\begin{equation}
  \label{eq:error_2}
  \Delta_n= 
  \frac{1}{2\pi i} \int_{C} d\zeta\; A_n(\zeta)\;
  \frac{1}{a^2 + \zeta^2} = 
  \frac{A_n(i\pi a)}{a}
\end{equation}
For large values of $a$ using the asymptotic expression for $A_n$ for
$\nu^2\gg \pi a$ and $a \gg \nu$ this expression simplifies to:
\begin{equation}
  \label{eq:B_large_tau}
  \Delta_n \approx 8 \nu 
  \exp\left(
    -\frac{\nu^2}{\pi a}
  \right)
\end{equation}

\begin{figure}
  \centering
  \includegraphics[angle=-90,width=\textwidth]{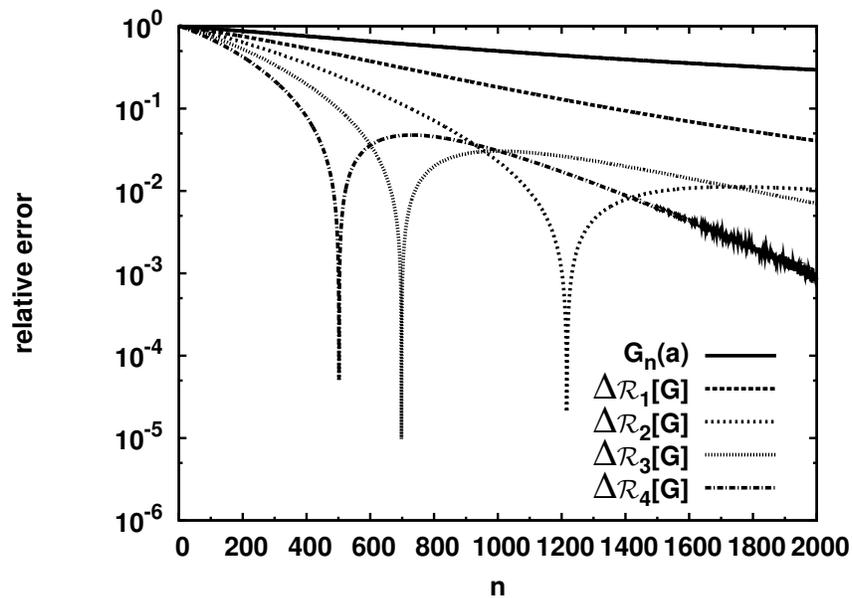}
  \caption{Relative error of the Richardson extrapolations of the 
  example 2, sum Eq.~\ref{eq:example_2}, for $a=1000$. Note the appearance
  of rounding errors in $\Delta\mathcal{R}_4[G](n)$.}
  \label{fig:Richardson}
\end{figure}

We contrast this behavior with other well known techniques for acceleration of
sums like Richardson extrapolation \cite{BenderOrzag}.  Using the
Euler-MacLaurin sum formula it can be shown that the partial sums behave
asymptotically as:
\begin{equation}
  \label{eq:partial_sum_asymptotic}
  G_n(a) = \frac{1}{a^2} + 2\sum_{k=1}^{n}\frac{1}{a^2+k^2} = 
  \frac{\pi\coth(\pi a)}{a}-\frac{2}{n}+\frac{1}{n^2}+
  \frac{2a^2-1}{3n^3}-\frac{a^2}{n^4}-\frac{6a^4-10a^2+1}{15n^5}+
  \ldots
  %O\left(\frac{1}{n^6}\right)
\end{equation}
The $N$th Richardson extrapolation for the partial sums at some $n$ is given
by:
\begin{equation}
  \label{eq:Richardson}
  \mathcal{R}_N[G] = 
  \sum_{k=0}^N \frac{(-1)^k (-1)^{k+N} (n+k)^N}{k!(N-k)!} G_{n+k}(a).
\end{equation}
The asymptotic $k^{-2}$ behavior of the sum, Eq.~\ref{eq:example_2}, is only
reached when $|k| \approx |a|$ and only then convergence of the partials sums
starts.  Applying Richardson extrapolation the terms in $1/n$ are successively
removed.  However because the coefficients of the $n^{-m}$ terms are
approaching $a^m$ the convergence is very slow and the asymptotic behavior
only sets in when $n\gg a$. This can be seen in the plot of the error in the
Richardson extrapolations of Eq.~\ref{eq:example_2},
Fig.~\ref{fig:Richardson}. The relative error is given by $\Delta
\mathcal{R}_N[G](a) = |(\mathcal{R}_N[G](a)-G(a))/G(a)|$.  It is clear that
Gaussian summation achieve the same accuracy with a substantially smaller
number of function evaluations. In fact the Richardson extrapolation shows
spurious convergence for partials sums with $N\sim{a}$.

%--- Distribution of zeros ----------------------------------------------------

\section{Distribution of zeros of the orthogonal polynomials}
The distribution of zeros of the denominators of the Pade approximants
$S_n(x)$ turns out to be quite unusual. According to
Eq.~\ref{eq:solution_S} the zeros of $S_N(x)$ for a given odd integer
number $N$ are determined by the equation:
\begin{equation}
  \label{eq:zeros}
  \cos(x)J_{n+5/2}(x)+\sin(x)Y_{n+5/2}(x) = 0.
\end{equation}
We are interested in the asymptotic behavior of the zeros at large
values of $n$.  Let $x>0$. The function is symmetric with respect to
$x\rightarrow-x$ so for each zero $x_0$ there is a corresponding zero
at $-x_0$. We define $\nu=n+5/2$ and $\tau=x/\nu$. The condition for a
zero of the denominator apart from non vanishing prefactors is:
\begin{equation}
  \label{eq:zeros_tau}
  \cos(\nu\tau)J_{\nu}(\nu\tau)+\sin(\nu\tau)Y_{\nu}(\nu\tau) = 0;
\end{equation}
In the limit $\nu\rightarrow\infty$ we have to differentiate between
two regimes, $\tau<1$ and $\tau>1$. For the first case, $\tau < 1$ the
leading asymptotic behavior of the Bessel function corresponding to
$x<\nu$ is given by \cite{Abramowitz,Gradstein}:
\begin{eqnarray}
  \label{eq:J_asymptotic}
  J_{\nu}(\nu\tau) &=& 
  \frac{1}{\sqrt{2\pi\nu}} 
  \frac{\exp(\nu(\sqrt{1-\tau^2}-{\rm arcosh}(1/\tau)))}
  {(1-\tau^2)^{1/4}}
  % \times \left[1+O\left(\frac{1}{\nu}\right)\right].
\end{eqnarray} and
\begin{eqnarray}
  \label{eq:Y_asymptotic}
  Y_{\nu}(\nu\tau) &=& 
  \frac{-2}{\sqrt{2\pi\nu}} 
  \frac{\exp(-\nu(\sqrt{1-\tau^2}-{\rm arcosh}(1/\tau)))}
  {(1-\tau^2)^{1/4}}
  % \times \left[1+O\left(\frac{1}{\nu}\right)\right].
\end{eqnarray}
The argument of the exponent of the second term is negative in the interval
$(0,1)$. The first term is exponentially suppressed compared to the first term.
The zeros of the denominator for $\tau\in(0, 1)$ are therefore determined
by the expression:
\begin{equation}
  \label{eq:zeros_tau_smaller_1}
  \sin(\nu\tau) = 0. 
\end{equation}
From this equation we conclude that $\nu\tau=\pi m$ with $m$ integer
of the same order as $\nu$ so it is useful to introduce $\sigma = m/\nu$, so
that in this regime $\tau = \pi \sigma$.

Next we consider the case $\tau>1$. The leading term in the asymptotic
expansion of the Bessel functions takes a different form:
\begin{eqnarray}
  \label{eq:J_asymptotic_tau_greater_1}
  J_{\nu}(\nu\tau) &=& \sqrt{\frac{2}{\pi\nu}}
  \cos(\nu(\sqrt{\tau^2-1}-\arccos(1/\tau))) 
  % + O\left(\frac{1}{\nu}\right)
  \\
  Y_{\nu}(\nu\tau) &=& \sqrt{\frac{2}{\pi\nu}}
  \sin(\nu(\sqrt{\tau^2-1}-\arccos(1/\tau))). 
  % + O\left(\frac{1}{\nu}\right)
\end{eqnarray}
For this case the leading asymptotic behavior of the zeros of the
denominator have to be determined by the following equation:
\begin{equation}
  \label{eq:zeros_tau_greater_1}
  \cos(\nu(\tau-\sqrt{\tau^2-1}+\arccos(1/\tau)) = 0.
\end{equation}
In this regime the zeros are solutions of 
$\tau-\sqrt{\tau^2-1}+\arccos(1/\tau)=\pi \sigma$.
To summarize the asymptotic behavior of the zeros is given by:
\begin{equation}
  \label{eq:zeros_g}
  \pi\sigma(\tau) =
  \begin{cases}
    \tau & \tau \le 1\\
    \tau-\sqrt{\tau^2-1}+\arccos(1/\tau)& \tau > 1.
  \end{cases}
\end{equation}
Note that $\sigma$ approaches $1/2$ in the limit $\tau\rightarrow\infty$.
This is the expected because we have considered the case $x>0$ only. In fact
the distribution of zeros is symmetric in the range $\sigma \in (-1/2,1/2)$
since the $S_n(x)$ depend $x^2$ so only half of the zeros lies above 0.  We
enumerate the positive zeros by $x_k$, $k\in\{1, 2, \ldots n\}$. The $n$
negative zeros are given by $-x_k$.  The asymptotic behavior of the zeros is
surprising: For large $n$ the first $[n/\pi]$ zeros are with high accuracy
given by $x=\pi, 2\pi, 3\pi, ...$. However for $k>[n/\pi]$ the distance
between two consecutive zeros become larger and eventually for very large
$x_k$ becomes very large and behaves like $x_k\sim \nu^2/{\pi(\nu-k)}$.  The
asymptotic behavior has been verified in numerical determination of the zeros.
This result is unusual because a singularity appears in the distribution of
zeros in the range of integration which is not present in the moment
generating function. The Gaussian summation achieves its accuracy by
evaluating $[n/\pi]$ points at the summation points and using the rest of the
points to explore the asymptotic behavior of the function which is being
summed. 

\section{Conclusions}
The general problem of summation of smooth functions is considered.
Summation formulas in analogy to Gaussian quadrature have been
derived.  The Pade approximants of the moment generating functions
have been used to derive the explicit solution of the recursion
relation. We have used this to give an explicit formula for the
corresponding orthogonal polynomials. The asymptotic distribution of
the zeros of these polynomials is derived.  We have shown that the
distribution of zeros shows an unusual behavior exhibiting a cusp. We
demonstrate that Gaussian summation needs a substantially smaller
number ($\sqrt{n}$ compared to $n$) of function evaluations for smooth
functions with a large scale than other common techniques. In
particular higher order Richardson extrapolation will lead to rounding
errors spoiling the convergence.

It would be interesting to investigate Gaussian
summation for different kind of linear functionals like
$\sum_{\nu\in\mathbb{Z}\backslash\{0\}}1/\nu^3\;f(1/\nu^2)$.
Summation schemes for summation over odd numbers (Matsubara
frequencies) relevant for semiclassical field theories in physics will
be discussed in a separate paper \cite{Monien}. Another interesting
point is the development of Kronrod schemes for summation which we are
currently working on.

\section{Acknowledgments}

I would like to acknowledge useful discussions with Stefan Engblom and
thank him in particular for writing his technical report
\cite{Engblom} on this interesting subject and providing me with
reference \cite{Milovanovic}.  During the initial stage of this work I
obtained useful comments from W. Gautschi which were greatly
appreciated.  I would like to thank David Heilmann, Simon Trebst and
Gang Li for carefully reading the manuscript. Finally I would like to
thank Andrew J. Millis for useful discussions and his remark ``If you
want a number then calculate.''. I would like to thank the {\em Aspen
  Center for Physics} where part of the work has been done.

\newpage

\bibliographystyle{amsplain}
\bibliography{gauss}

\providecommand{\bysame}{\leavevmode\hbox to3em{\hrulefill}\thinspace}
\providecommand{\MR}{\relax\ifhmode\unskip\space\fi MR }
\providecommand{\MRhref}[2]{%
  \href{http://www.ams.org/mathscinet-getitem?mr=#1}{#2}
}
\providecommand{\href}[2]{#2}
\begin{thebibliography}{10}

\bibitem{Baker}
George A.~JR. Baker, \emph{Essentials of {P}ade {A}pproximants}, Academic
  Press, NY, 1975.

\bibitem{BenderOrzag}
Carl~M. Bender and Steven~A. Orzag, \emph{Advanced {M}athematical {M}ethods for
  {S}cientists and {E}ngineers}, McGraw-Hill, 1984.

\bibitem{Engblom}
Stefan Engblom, \emph{Gaussian quadratures with respect to discrete measures},
  Uppsala University, Technical Report 2006-007 (2006).

\bibitem{Gautschi}
Walter Gautschi, \emph{Orthogonal {P}olynomials}, Oxford University Press,
  2004.

\bibitem{GautschiHL}
\bysame, \emph{The {H}ardy-{L}ittlewood {F}unction: {A}n {E}xercise in {S}lowly
  {C}onvergent {S}eries}, J. Comput. Appl. Math. \textbf{179} (2005), 249--254.

\bibitem{Golub}
G.~H. Golub and J.~H. Welsh, \emph{Calculation of {G}auss quadrature rules},
  Math. Comp. \textbf{23} (1969), 221--230.

\bibitem{Gradstein}
I.~S. Gradshteyn and I.~M. Ryzhik, \emph{Table of {I}ntegrals, {S}eries and
  {P}roducts}, Academic Press, 1994.

\bibitem{Homeier}
Herbert~H.~H. Homeier, \emph{Scalar {L}evin-type sequence transformations}, J.
  Comput. Appl. Math. \textbf{122} (2000), 81--147.

\bibitem{Levin}
D.~Levin, \emph{Methods for accelerating convergence of infinite series and
  integrals}, Ph.D. thesis, Tel Aviv University, 1975.

\bibitem{Abramowitz}
{M. Abramowitz and Irene A. Stegun} (ed.), \emph{Handbook of {M}athematical
  {F}unctions}, Dover, 1965.

\bibitem{Milovanovic}
{Gradimir~V.} Milovanovi{\'c} and Aleksandar Cvetkovi{\'c}, \emph{Convergence
  of {G}aussian {Q}uadrature {R}ules for {A}pproximation of {C}ertain {S}ums},
  Eastern J. on Approximation \textbf{10} (2004), 171--187.

\bibitem{Monien}
Hartmut Monien, \emph{Semiclassical mean field theories}, unpublished.

\bibitem{Nikiforov}
Arnold~F. Nikiforov and Vasilii~B. Uvarov, \emph{Special {F}unctions of
  {M}athematical {P}hysics}, Birkh\"auser Verlag Basel, 1978.

\bibitem{NumericalRecipes}
William~H. Press, Saul~A. Teukolsky, William~T. Vetterling, and Brian~P.
  Flannery, \emph{Numerical {R}ecipes in {C}++}, Cambridge University Press,
  2002.

\end{thebibliography}

\begin{table}
\begin{tabular}{|c|c|c|c|c|c|c|}
\hline
n & $\Delta H(1)$ & $\Delta H(5)$ & $\Delta H(10)$ & $\Delta H(20)$ & $\Delta H(40)$ & $\Delta H(100)$ \\  \hline
2 & $8.73\cdot 10^{-9}$ & $1.58\cdot 10^{-2}$ & $1.9\cdot 10^{0}$ & $4.47\cdot 10^{-1}$ & $5.06\cdot 10^{-1}$ & $5.61\cdot 10^{-1}$ \\ 
3 & - & $2.53\cdot 10^{-6}$ & $1.02\cdot 10^{-2}$ & $9.55\cdot 10^{-1}$ & $2.29\cdot 10^{-1}$ & $3.29\cdot 10^{0}$ \\ 
4 & - & $3.66\cdot 10^{-11}$ & $3.3\cdot 10^{-6}$ & $1.67\cdot 10^{-2}$ & $1.28\cdot 10^{0}$ & $3.29\cdot 10^{0}$ \\ 
5 & - & - & $1.47\cdot 10^{-10}$ & $2.29\cdot 10^{-5}$ & $2.48\cdot 10^{-1}$ & $4.51\cdot 10^{-1}$ \\ 
6 & - & - & - & $5.19\cdot 10^{-9}$ & $3.04\cdot 10^{-3}$ & $2.53\cdot 10^{0}$ \\ 
7 & - & - & - & $2.85\cdot 10^{-13}$ & $5.89\cdot 10^{-6}$ & $2.77\cdot 10^{0}$ \\ 
8 & - & - & - & - & $2.8\cdot 10^{-9}$ & $1.09\cdot 10^{0}$ \\ 
9 & - & - & - & - & $4.19\cdot 10^{-13}$ & $4.46\cdot 10^{-2}$ \\ 
10 & - & - & - & - & - & $3.87\cdot 10^{-4}$ \\ 
11 & - & - & - & - & - & $1.02\cdot 10^{-6}$ \\ 
12 & - & - & - & - & - & $1.01\cdot 10^{-9}$ \\ 
13 & - & - & - & - & - & $4.07\cdot 10^{-13}$ \\ 
14 & - & - & - & - & - & $2.51\cdot 10^{-14}$ \\ 
15 & - & - & - & - & - & - \\ 
\hline
\end{tabular}

  \caption{Results for the relative error in evaluation of the 
    Hardy-Littlewood sum.}
  \label{tab:hardy_littlewood_sum_error}
\end{table}

\begin{figure}
  \centering
  \includegraphics[angle=-90,width=\textwidth]{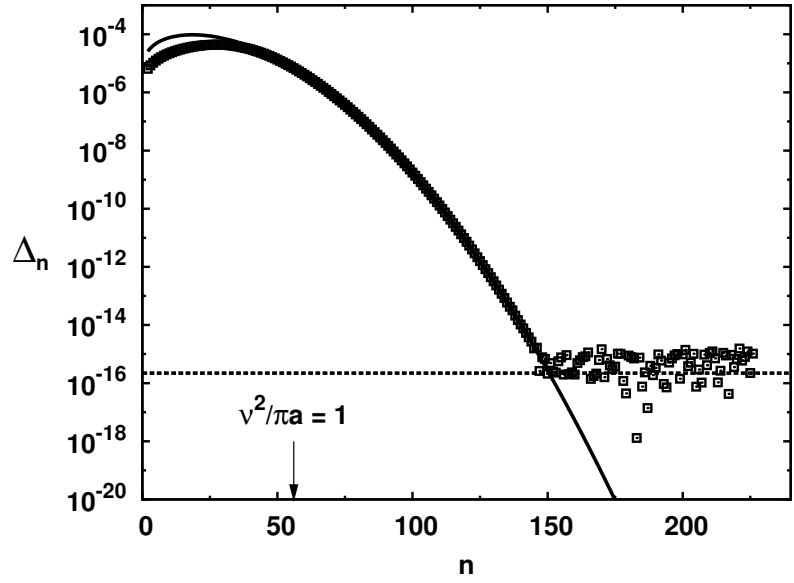}
  \caption{Relative error of the Gaussian summations of the 
    example 2, sum Eq.~\ref{eq:example_2}, for $a=1000$. The solid line is
    the error estimate given in the text ($\sim 8n\exp(-4n^2/\pi{a})$) and the
    dotted line is giving the machine $\epsilon$.  }
  \label{fig:convergence_as_function_of_n}
\end{figure}

%% \begin{figure}
%%   \centering
%%   \includegraphics[width=\textwidth]{extrapolation}
%%   \caption{Number of function evaluations needed to
%%     achieve machine precision for example 2, sum Eq.~\ref{eq:example_2}, for
%%     variing $a$ (solid line) and the corresponding numerical result of fitting
%%     the logarithm of the error with a quadratic function (squares).}
%%   \label{fig:convergence_as_function_of_a}
%% \end{figure}

\begin{figure}
  \centering
  \includegraphics[angle=-90,width=\textwidth]{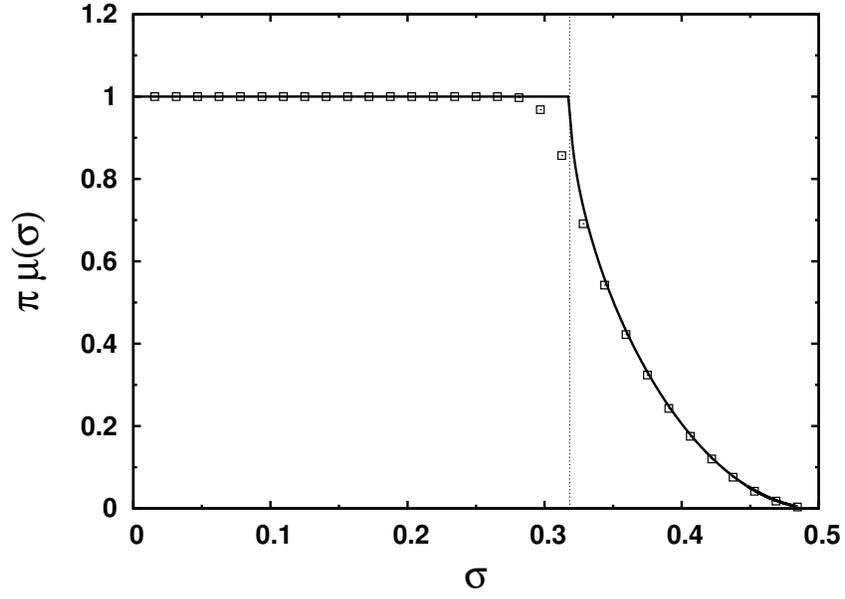}
  \caption{Asymptotic density of zeros (
    $\mu(\sigma)=\frac{d\tau}{d\sigma}(\sigma)$). The dots give the
    approximate density of zeros $1/(x_{n+1}-x_{n})$ for $N=32$ plotted as a
    function of $\sigma=n/N$. The dotted line indicates $\sigma=1/\pi$.  }
  \label{fig:density_of_zeros}
\end{figure}

%% \begin{figure}
%%   \centering
%%   \includegraphics[width=\textwidth]{convergence_matsubara}
%%   \caption{Relative error in the evaluation of the sum Eq.~\ref{eq:action_sum}
%%     for $x=1000$.  The dashed line gives the machine $\epsilon$ for
%%     comparison.}
%%   \label{fig:matsubara_convergence}
%% \end{figure}

\end{document}